# Infinite Set of Resonances in the Linear Damped Oscillator Subject to Harmonic Forcing with Non-standard Frequency Modulation

Oleg V. Gendelman[1], Luis M. Baldelomar Pinto[2*],
Lawrence A. Bergman[3], Alexander F. Vakakis[2]

[1] Faculty of Mechanical Engineering, Technion – Israel Institute of Technology Haifa, Israel
[2] University of Illinois Urbana-Champaign, Mechanical Science and Engineering, Urbana, USA
[3] University of Illinois Urbana-Champaign, Aerospace Engineering, Urbana, USA

*Corresponding author, luismb2@illinois.edu

## Abstract

It is shown that harmonic signals incorporating a type of weak non-standard frequency modulation (wNSFM) have interesting spectral properties, namely, time-dependent bandwidths that become increasingly broader with increasing time. As such, they represent a class of signals with frequency-time coupling in their spectra. Applying this type of wNSFM harmonic excitation to the classical single-degree-of-freedom (SDOF) linear, time-invariant (LTI) damped/undamped oscillator yields interesting resonance phenomena. Specifically, the weakly damped oscillator forced by a wNSFM harmonic excitation with a central frequency not necessarily equal to its natural frequency, exhibits two transient resonance captures involving two distinct harmonics with relatively high amplitudes over finite time intervals, while the long-term response decays as $\sim t^{-1/2}$ as $t \to \infty$. Moreover, this set of resonance captures is persistent to changes in the central frequency of the force modulation, even if it lies well apart from the natural frequency of the oscillator. Considering the corresponding undamped oscillator, it possesses two types of resonances, referred to as simple and non-simple resonances. Simple resonances correspond to finite-amplitude steady-state responses caused by two sustained resonance captures, in the form of two superimposed distinct quasi-periodic responses, which, however, are "activated" at different time instants. Again, this set of sustained resonances is persistent to changes in the frequency modulation, in the sense that even if a given resonance is eliminated by parameter variation, another one is realized in its neighborhood. The necessary and sufficient conditions for non-simple resonances are given in the form of a theorem which predicts the co-existence of as many as four possible resonant harmonics and specifies special conditions that the phases of these resonant harmonics must satisfy for constructive interference; the resulting undamped non-simple resonance grows unboundedly as $\sim t^{1/2}$ as $t \to \infty$, in contrast to the classical linear-in-time ($\sim t$ as $t \to \infty$) resonance growth of the linear resonator with unmodulated harmonic excitation. Similar to the simple resonances, the set of non-simple resonances is persistent to changes in the central frequency and the other parameters of the wNSFM. Our results reveal interesting infinite sets of resonances in linear SDOF resonators under frequency-modulated excitations.

## 1. Introduction

It is typical in dynamics to study the responses of systems subject to time-periodic forces, which can be periodic, quasi-periodic or even chaotic, through analytical and numerical techniques **[1-5]**. Intuitively then, one can ascertain that modulating in some sense the time-periodicity in the forcing should result in more complicated and rather unexpected dynamics. In a recent work **[6]** it was shown that applying harmonic excitation with *modulated amplitude* can enable nonlinear hysteresis loop "tracking" in a Duffing oscillator; also, it can eliminate the well-known phenomenon of *bi-stability*, i.e., the uncertainty related to the attraction of the steady-state dynamics by different co-existing stable attractors depending on the initial conditions. No similar effects, however, were observed when *frequency modulation* of the same harmonic excitation was considered, a result which motivated the authors to conduct the present study.

Frequency modulation was developed mainly in telecommunications and electronics for transmission of information through radio waves, whereby the instantaneous frequency of a "carrier wave" is varied in time, e.g., periodically, in proportion to a certain perturbation amplitude **[7]**. Typically, the frequency modulated signal has a form of the type,

$$g_1(t) \equiv \cos(\omega_0 t + \varepsilon\alpha \sin\Omega t) \quad (1)$$

where $\omega_0$ is the central frequency of the carrier signal, and $\varepsilon\alpha$, $\Omega$, are the amplitude and frequency parameters, respectively, that characterize the frequency modulation. Assuming *weak frequency modulation*, we consider $0 < \varepsilon \ll 1$ to be a small parameter and all other parameters to be of $\mathcal{O}(1)$. This is the *standard form of weak frequency modulation*, having an instantaneous phase $\phi_1(t) = \omega_0 t + \varepsilon\alpha \sin(\Omega t)$, and an instantaneous frequency $\dot{\phi}_1(t) = \omega_0 + \varepsilon\alpha\Omega \cos\Omega t$; this represents a small time-periodic perturbation of the center frequency. Note that since (1) is, in general, a quasi-periodic time signal, its frequency spectrum consists of an infinite number of "side bands" above and below the central frequency $\omega_0$, which, however, have increasingly decaying amplitudes away from $\omega_0$, so higher or lower harmonics may be neglected for all practical purposes. In addition, the responses of linear resonators forced by weakly modulated harmonic excitations of the form (1) can be analytically studied through asymptotic techniques **[4]**. Some additional works considering this type of frequency modulation are **[8-12]**.

In this work a non-standard form of frequency modulation, is introduced by weakly modulating directly the central frequency as follows:

$$g_2(t) \equiv \cos[(\omega_0 + \varepsilon\alpha \sin\varepsilon t)\, t], \quad 0 < \varepsilon \ll 1 \quad (2)$$

We denote this as *weak non-standard frequency modulation (wNSFM)* where the characterization "weak" refers to both the amplitude and frequency of the perturbation of the central frequency $\omega_0$. For the signal (2) the instantaneous phase is $\phi_2(t) = (\omega_0 + \varepsilon\alpha \sin\varepsilon t)\, t$, and the instantaneous frequency is $\dot{\phi}_2(t) = \omega_0 + \varepsilon\alpha \sin\varepsilon t + \varepsilon^2 t\alpha \cos\varepsilon t$; moreover, one can define also an effective (average) frequency as $\phi_2(t)/t = \omega_0 + \varepsilon\alpha \sin\varepsilon t$. Conceptually, the instantaneous and effective (average) frequencies could be regarded as analogous to the well-known phase and group velocities defined for traveling waves in acoustics **[13]**. Interestingly, the instantaneous frequency of (2) is oscillatory with an envelope that increases linearly in time, so it could be regarded as being in a state of linear *phase resonance*.

A special note concerns the results reported in **[8]**, and their relation to the present wNSFM signal (2). Note that (2) can be written as, $g_2(t) \equiv \cos[(\omega_0 + \varepsilon\alpha \sin\varepsilon t)\, t] = \cos(\omega_0 t + I(t)\sin\varepsilon t)$, where $I(t) = \varepsilon\alpha t$, which can be regarded as a *time-varying modulation index* in the notation of **[8]** compared to the *constant modulation index* $I = \varepsilon\alpha$ in the standard frequency modulation (1) – and equation (1) in **[8]**. Note that Chowning **[8]** did contemplate the possibility of a time-varying modulation index – he called the resulting spectrum, "dynamic spectrum" – and made the following statement in his work: "…The complexity of the spectrum is related to the modulation index in such a way that as the index increases, the bandwidth of the spectrum also increases (see Fig. 5). If, then, the modulation index were made to be a function of time, the evolution of the band width of the spectrum could be generally described by the shape of the function…" So, Chowning contemplated the possibility of making the modulation index time-dependent (as is performed in the present work), but this point was not pursued rigorously in his work (nor in any other, to the authors' best knowledge) except for some numerical examples of "dynamic spectra" and relevant comments given in **[8]**. Hence, the present analytical study of the time-dependent bandwidth of the signal (2), and even more importantly, the dynamic analysis of the classical harmonic oscillator subject to the forcing (2) with wNSFM, represent new contributions to the literature.

Another observation regarding the wNSFM signal (2) stems from its comparison with a well-known signal such as the linear "chirp", $c(t) = A\, cos\left(\omega_0 t + \frac{1}{2}kt^2\right)$. Note that $c(t)$ has *a single fast frequency* $\omega_0$, but it also produces time-frequency coupling in the Fourier spectrum. A simple comparison with the wNSFM signal (2), $g_2(t) \equiv \cos[(\omega_0 + \varepsilon\alpha \sin\varepsilon t)\, t],\ 0 < \varepsilon \ll 1$, shows that these are completely different signals; moreover, whereas the chirp results in a *monotonically increasing (or decreasing)* phase, the wNSFM yields an *oscillating phase with a continuously increasing envelope*. As such, in (2) an *additional slow frequency* $\varepsilon$ is "embedded" into the signal; combined, these *two frequencies* (the fast central frequency $\omega_0$ and the slow modulation frequency $\varepsilon$) generate an infinity of harmonics in the spectrum of the signal, which, as shown in this work, produce new types of combination resonance phenomena even in the simplest of dynamical systems, namely the classical linear, time-invariant oscillator; such resonance phenomena cannot be achieved with the chirp, which produces only passage through resonance **[9,10]**. To our best knowledge, the realization of combination resonances in the classical harmonic oscillator forced by a single frequency-modulated excitation is introduced for the first time in the field of dynamics and vibrations.

Lastly, the signal (2) is special for an additional reason. Regular amplitude modulation – considered, e.g., in **[6]**, or common phase modulation, e.g., eq. (1), are, in essence, narrowband processes, so broadband excitations are usually associated with applied excitations over a broad frequency range, white noise being the most common example. In contrast, the wNSFM excitation (2) strongly differs from these well-known patterns, since a broad frequency range is achieved over a substantial time span and in a well-controlled manner. As a result, a controlled and predictable broadband frequency-time coupling effect can be achieved, with the resulting band being decided by the selected window of time integration. This is the reason why the resulting response patterns discussed below are not achievable by other standard forcing conditions, even in the simplest case of the classical single-degree-of-freedom (SDOF) linear oscillator.

In the next section it is shown that the signal (2) with wNSFM has interesting spectral properties, the most important being that its Fourier spectrum is time-dependent and becomes increasingly more broadband with increasing time. In that sense, it can be regarded as a highly non-stationary signal, as it has a continuously expanding broadband frequency spectrum with time, as an increasingly larger number of higher harmonics is "activated" as time increases. Hence, signals of the form (2) appear to possess the interesting feature of strong *frequency-time coupling*, which is not typically encountered in dynamics, being more of a feature in the realm of acoustics, e.g., in traveling waves where frequency-wavenumber-time coupling exists. Then, the dynamics of the classical SDOF, linear time invariant (LTI) oscillator subject to harmonic excitations with wNSFM (2) will be studied, and new types of persistent transient and sustained multi-harmonic resonance phenomena will be analyzed in detail and validated numerically.

## 2. Harmonic signals with weak non-standard frequency modulation (wNSFM)

To study the spectral content of signal (2), the following Jacobi-Anger expansions are used whereby exponentials of trigonometric functions are expressed in terms of their harmonics **[14]**,

$$e^{jz\cos\theta} = \sum_{n=-\infty}^{\infty} j^n J_n(z)\, e^{jn\theta}, \quad e^{jz\sin\theta} = \sum_{n=-\infty}^{\infty} J_n(z)\, e^{jn\theta} \tag{3}$$

where $J_n(z)$ is the $n$ −th order Bessel function of the first kind, $z$ is a complex argument, and $j = (-1)^{1/2}$. From (3), after some straightforward manipulations one can obtain the following harmonic expansion of signal (2):

$$\begin{gathered} g_2(t) \equiv \cos[(\omega_0 + \varepsilon\alpha \sin\varepsilon t)\, t] = \\ \frac{1}{2}\left[e^{j\omega_0 t}\, e^{j\varepsilon\alpha t \sin\varepsilon t} + e^{-j\omega_0 t}\, e^{-j\varepsilon\alpha t \sin\varepsilon t}\right] = \\ \sum_{n=-\infty}^{\infty} J_n(\varepsilon\alpha t)\, cos\,(\omega_0 + \varepsilon n)t \end{gathered} \tag{4}$$

Then, the mechanism of evolution in time of the Fourier spectrum of (2) becomes clear. Indeed, for small values of their arguments, Bessel functions of the first kind obey the following asymptotic behavior, $n\, J_n(x) \sim \frac{1}{n!}\left(\frac{x}{2}\right)^n$, $0 < x \ll 1$. It follows that since $0 < \varepsilon \ll 1$, an increasing number of higher harmonics in (4) (characterized by large values of $n$) is "activated", i.e., their amplitudes achieve $\mathcal{O}(1)$ values, over an increasingly extended time span. In fact, to get a leading-order analytic estimation for the activation of the $n$ −th harmonic, one can invoke Stirling's formula, $n! \sim (n/e)^n$ valid for large $n$, to predict that the amplitude of the $n$ −th harmonic becomes of $\mathcal{O}(1)$ at the characteristic time $\tau_n \sim n/\varepsilon\alpha$, which, in turn, indicates that the "frequency cutoff" for the $n$ −th harmonic then can be approximated by $\omega_{cutoff,n} \sim \varepsilon n$, or taking into account the previous estimate for the activation time, $\omega_{cutoff,n} \sim a\varepsilon^2\tau_n$. Hence, the wNSFM signal (2) possesses a *time-dependent bandwidth*, in the sense that its Fourier spectrum continuously expands in frequency with progressing time, and increasingly more delayed higher harmonics are activated as time increases. This is in contrast to the standard frequency modulation (1) whose spectral properties are time-invariant.

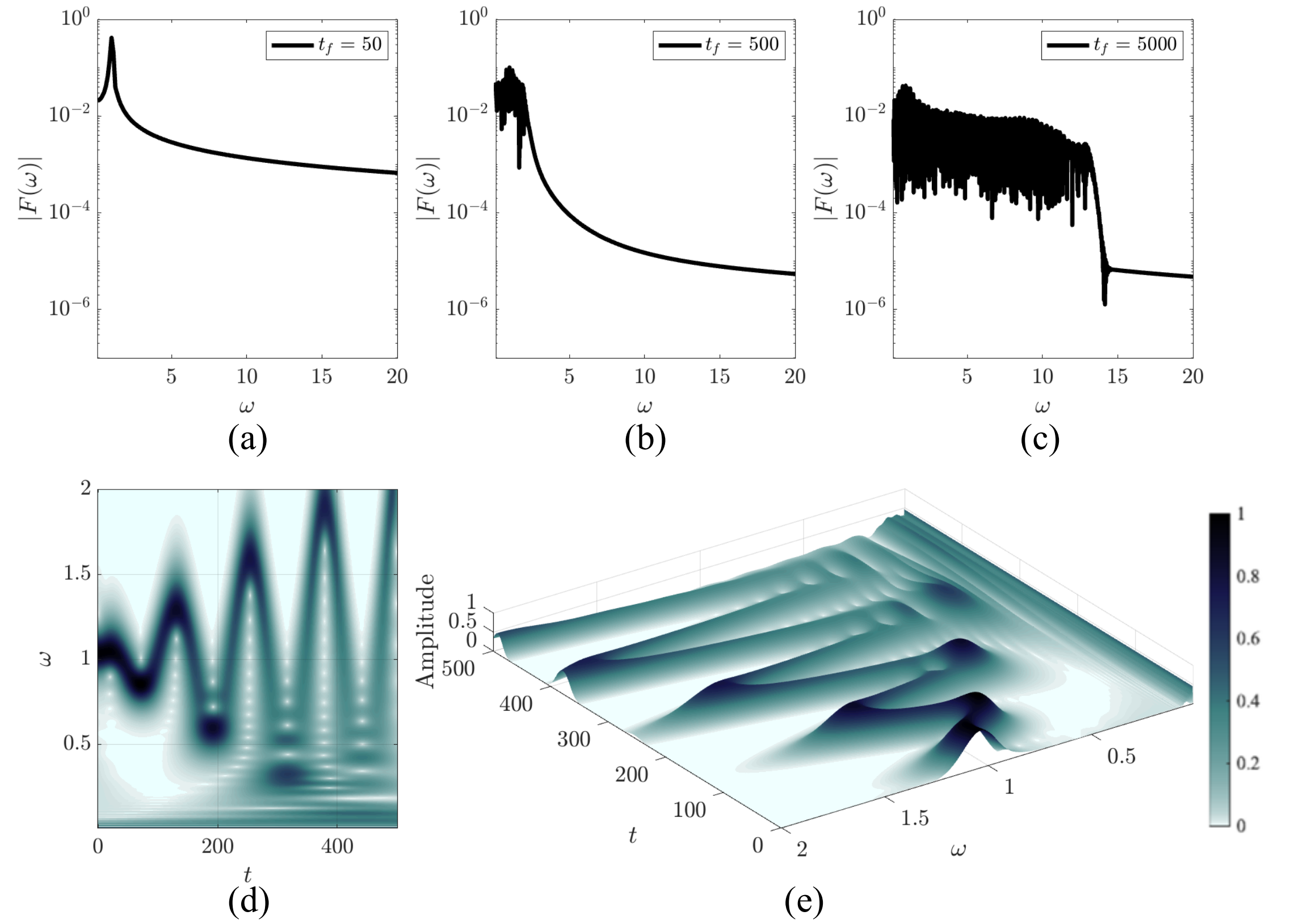

Figure 1. Spectral study of signal (4) with wNSFM for $\omega_0 = 1, \varepsilon = 0.05$ and $\alpha = 1$: Fourier spectrum for integration time of, (a) $\tau_f = 50$, (b), $\tau_f = 500$, and (c) $\tau_f = 5{,}000$; (d) contour plot of the wavelet transform spectrum, and (e) wavelet transform spectrum.

In Fig. 1, the time-varying spectral content of signal (4) is depicted for a specific set of system parameters, and the delayed activation of higher harmonic components is clearly discerned, as well as the time-frequency coupling in the spectral composition of this signal. One notes the continuous generation of higher harmonics with increasing integration time, which is a characteristic feature of the wNSFM (4) and differentiates it clearly from the signal (1) with conventional frequency modulation which has a definitive frequency spectrum with set side bands left and right of the center frequency. In the following sections the responses of damped or undamped linear oscillators subject to excitations with non-standard frequency modulation will be studied, and new types of non-conventional transient and sustained resonance phenomena (resonance captures) are revealed. Moreover, it is shown that these resonances persist with changes in the modulation parameters, and even with variations of the central frequency $\omega_0$.

## 3. Damped linear resonator subject to harmonic forcing with wNSFM

First, the response of the classical LTI SDOF oscillator subject to a harmonic excitation with wNSFM is studied for zero initial conditions; this does not restrict the generality of the analysis given that, due to the linearity of the problem, the homogeneous solution is of standard form. Hence, the following linear system is examined,

$$\ddot{x}(t) + \lambda\dot{x}(t) + x(t) = \cos[(\omega_0 + \varepsilon\alpha\,\sin\varepsilon t)\,t] = \sum_{n=-\infty}^{\infty} J_n(\varepsilon\alpha t)\, cos\,(\omega_0 + \varepsilon n)t,\;\; x(0) = \dot{x}(0) = 0,\;\; 0 < \varepsilon \ll 1 \tag{5}$$

with $0 < \lambda < 2$ for an underdamped oscillator with normalized natural frequency equal to unity (again, without loss of generality). Then, the forced solution can be expressed in closed form through the convolution integral,

$$x(t) = \frac{exp(\lambda t/2)}{\beta}\int_0^t exp(-\lambda\tau/2)\, sin\beta(t-\tau)\cos[(\omega_0 + \varepsilon\alpha\,\sin\varepsilon\tau)\,\tau]\, d\tau = \frac{exp(\lambda t/2)}{\beta}\sum_{n=-\infty}^{\infty}\int_0^t exp(-\lambda\tau/2)\, sin\beta(t-\tau)\, J_n(\varepsilon\alpha\tau)\cos(\omega_0 + \varepsilon n)\tau\, d\tau \tag{6}$$

where $\beta = \sqrt{1 - \lambda^2/4}$. To assess the asymptotic behavior of the solution as $t \to \infty$, one considers the asymptotic expressions of the Bessel functions of the first kind for large arguments and use integration by parts to obtain the following asymptotic expressions,

$$\int_0^t exp\left[\left(\frac{\lambda}{2} + j\omega\right)\tau\right] J_n(\varepsilon\alpha\tau)\, d\tau \sim \frac{exp\left(\frac{\lambda t}{2}\right)}{\sqrt{2\pi\varepsilon\alpha t}}\left\{\frac{exp[j(\varepsilon\alpha t - \varphi_n)]}{\frac{\lambda}{2} + j(\omega + \varepsilon\alpha)} + \frac{exp[-j(\varepsilon\alpha t - \varphi_n)]}{\frac{\lambda}{2} + j(\omega - \varepsilon\alpha)} + O\left(\frac{1}{t}\right)\right\},\;\; t \to \infty \tag{7}$$

where $\varphi_n = \frac{\pi n}{2} + \frac{\pi}{4}$.

Considering the asymptotic estimate (7), one concludes that each separate harmonic component of the infinite series expansion (6) decays as $\sim\sqrt{\frac{1}{t}}$ when $t \to \infty$. To claim that the complete infinite series (6) decays uniformly in the same way, one needs to demonstrate that the sums of coefficients of the inverse powers of $t$ do not diverge. This fact can be demonstrated by considering the asymptotic behavior of the single harmonic in the expansion (3) for large times and large frequencies. Generically it suffices to examine the asymptotic behavior of the solution of the following, somewhat simpler dynamical system,

$$\ddot{x}(t) + \lambda\dot{x}(t) + x(t) = \frac{exp(j\Omega t)}{\sqrt{t}},\;\; x(0) = \dot{x}(0) = 0 \tag{8}$$

for $t \gg 1$ and $\Omega \gg 1$. Indeed, the solution of (8) satisfies the following asymptotic behavior:

$$x(t) \sim exp(j\Omega t)\left[\frac{1}{\sqrt{t}(-\Omega^2 + j\lambda\Omega + 1)} + \mathcal{O}\left(\frac{1}{\sqrt{t}\Omega^4}\right)\right], \quad t, \Omega \gg 1 \tag{9}$$

Hence, comparing to the forced response (6) and letting $\Omega = \omega_0 + \varepsilon n$ one concludes that the infinite series (6) will behave as $(\sum_{n=-\infty}^{\infty} 1/n^2)/\sqrt{t}$, i.e., it will absolutely converge to zero as $t \to \infty$. Hence, the response is absolutely convergent, despite the fact that the Jacobi-Anger series do not possess this property. The conclusion on convergence is valid also for the system without damping, provided that the denominator in (9) does not vanish; this latter restriction isolates the resonant regimes.

From relations (8) and (9), one concludes that the steady-state solution (6) for the damped system will always decay to zero. Moreover, given that the homogeneous solution of (5) exponentially decays to zero, one proves that *the steady-state solution of the damped oscillator (5) always decays to zero irrespective of the parameters of the modulated excitation and the initial conditions*. A numerical demonstration of the decaying steady-state solution of system (5) is given in Fig. 2a for a resonant center frequency $\omega_0 = 1$. The equation of motion was integrated numerically using the MATLAB® ode45 solver. The solver employs an adaptive explicit Runge–Kutta method of order 4 based on the Dormand–Prince pair **[15]**. The relative and absolute error tolerances were set to $10^{-10}$ and $10^{-12}$, respectively, and the maximum integration time step was restricted to $10^{-3}$. The solver adaptively adjusted the time step while satisfying these prescribed error tolerances.

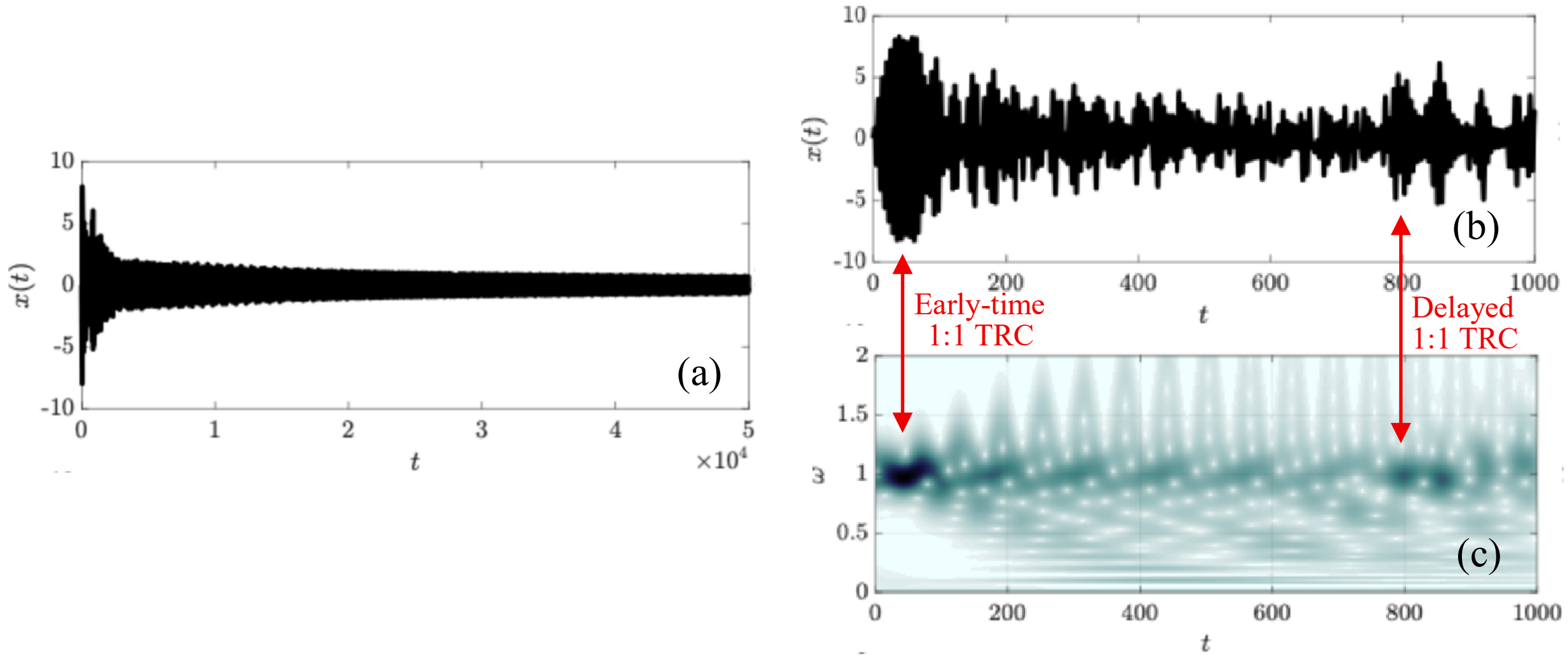

Figure 2. Forced response of (5) for $\varepsilon = 0.1$, $\alpha = 0.25$ and $\lambda = 0.1$ and resonant central frequency $\omega_0 = 1$: (a) Long-time solution showing the decaying steady state, (b) early-time solution, and (c) its corresponding wavelet transform spectrum, revealing two distinct 1:1 transient resonance captures (1:1 TRCs).

Note that in the absence of the weak modulation, and for relatively small viscous damping coefficient $\lambda$, the steady-state response would gain a large-amplitude resonance response *only* when the central frequency of the excitation satisfies the classical resonance condition $\omega_0 \approx 1$.

The weak modulation (4) prevents the formation of this type of classical resonance at steady state, but as discussed below, it produces new types of persistent transient resonance captures involving transient resonances of distinct harmonic (frequency) components only over finite time ranges.

From the early-time time series and corresponding wavelet spectrum of Figs. 2b,c one deduces two localized "bursts" at different time intervals, namely an initial high-amplitude burst of the response immediately following the application of the modulated harmonic excitation, followed by a second delayed burst at later time interval. Considering the governing equation of motion (5) it is evident that these two bursts correspond to transient resonance (or resonance captures) of two harmonics of the modulated force satisfying $\omega_0 + \varepsilon n = \pm 1$, yielding the two harmonics $n = (\pm 1 - \omega_0)/\varepsilon$, which for the specific parameters used in the simulations of Fig. 2 are two integer harmonics with $n = 0$, and $n = -20$. Hence, provided that the resulting values of $n$ in the previous resonance conditions are integers, a straightforward analysis predicts the existence of two distinct 1:1 transient resonance captures (1:1 TRCs) in the forced response of (5), both of which are due to their coincidence with the natural frequency of the oscillator, but only over a finite time interval, after which there occurs escape from resonance capture **[2]**. Moreover, the different activation of each of these harmonics in time (see discussion in the previous section), can be helpful in predicting the time of occurrence of each resonance following the application of the force; whereas the 20-th harmonic is activated with capture: Indeed, the zero-th harmonic is activated right from the start, so the first resonance capture occurs immediately approximate time delay $\tau_{20} \sim 20/\varepsilon\alpha = 800$, which agrees with the numerical simulations of Figs. 2b,c. It needs to be stated though that, *whereas these transient high-amplitude response bursts occur at the relatively early regime of the forced dynamics, still, the steady-state response of system (5) always decays to zero as* $t \to \infty$.

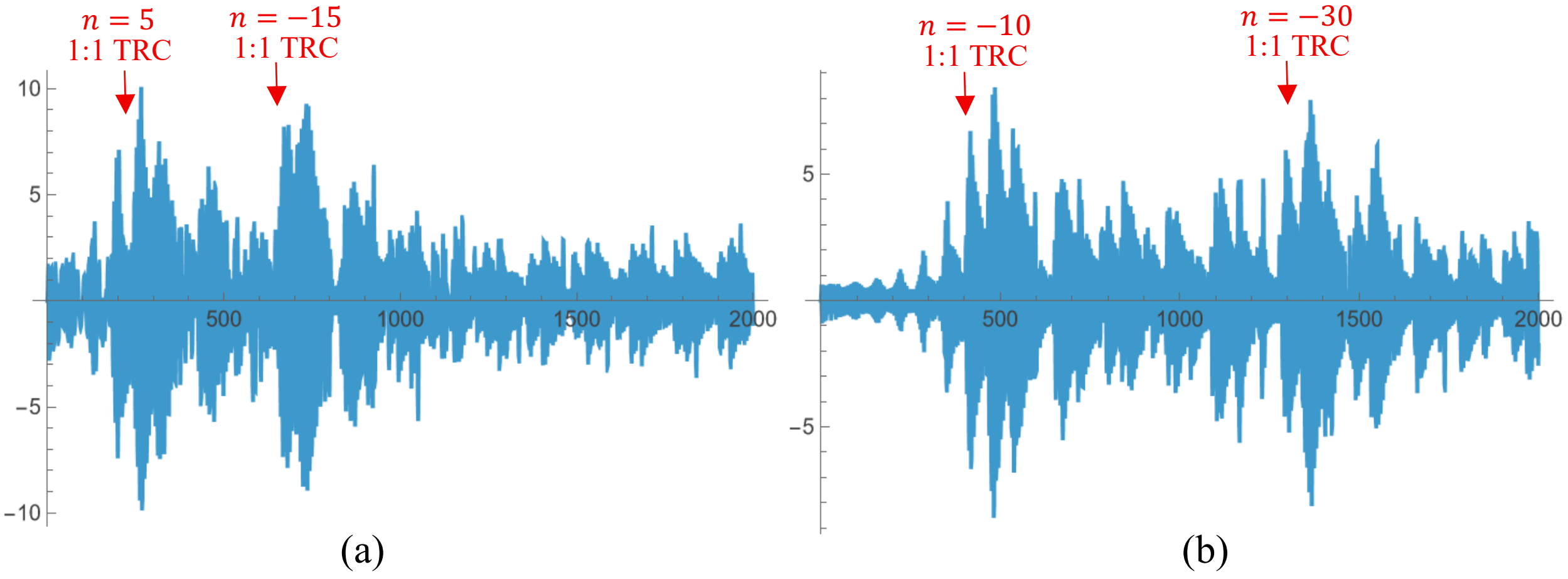

Figure 3. Persistent 1:1 transient resonance captures in (5) for $\varepsilon = 0.1$, $\alpha = 0.25$, $\lambda = 0.05$, and off-resonant center frequencies: (a) $\omega_0 = 0.5$, and (b) $\omega_0 = 2.0$.

Lastly, it is of interest to note that the two distinct 1:1 TRCs are *persistent* to changes in system parameters, that is, even when the center frequency of the wNSFM is itself *off-resonance*, i.e., well apart from the natural frequency of the oscillator, $\omega_0 \neq 1$. This can be clearly inferred from the 1:1 resonance condition, $n = (\pm 1 - \omega_0)/\varepsilon$, which, provided that $n$ is integer, can be satisfied for arbitrary off-resonance values of the center frequency, although then both 1:1 TRCs are delayed. This is shown in Fig. 3, where two different excitations with wNSFM and off-resonant center frequencies are considered, exhibiting both the predicted 1:1 TRCs even for center frequencies below (Fig. 3a) and above (Fig. 3b) the natural frequency of the oscillator equaling unity.

## 4. Undamped linear resonator subject to harmonic forcing with wNSFM

The steady-state response of the undamped system is now studied, which is governed by the differential equation:

$$\ddot{x}(t) + x(t) = \cos[(\omega_0 + \varepsilon\alpha \sin\varepsilon t)\, t] = \sum_{n=-\infty}^{\infty} J_n(\varepsilon\alpha t)\, cos\,(\omega_0 + \varepsilon n)t,\;\; x(0) = \dot{x}(0) = 0,\;\; 0 < \varepsilon \ll 1 \tag{10}$$

Since this is a linear dynamical system, its response is expressed in closed form by the convolution integral:

$$\begin{aligned} x(t) &= \int_0^t sin(t-\tau) \cos[(\omega_0 + \varepsilon\alpha \sin\varepsilon\tau)\, \tau]\, d\tau = \\ &\sum_{n=-\infty}^{\infty} \int_0^t sin(t-\tau)\, J_n(\varepsilon\alpha\tau) \cos(\omega_0 + \varepsilon n)\tau\, d\tau = \\ &sint \sum_{n=-\infty}^{\infty} \int_0^t J_n(\varepsilon\alpha\tau)\, cos\tau \cos(\omega_0 + \varepsilon n)\tau\, d\tau - \\ &cost \sum_{n=-\infty}^{\infty} \int_0^t J_n(\varepsilon\alpha\tau)\, sin\tau \cos(\omega_0 + \varepsilon n)\tau\, d\tau \end{aligned} \tag{11}$$

It is of interest to examine the long-term behavior of the solutions of (10), and the conditions under which resonance can occur. It will be shown that, depending on the parameters of the non-standard frequency modulation there are two cases of resonances that can be realized, which will be referred to as *simple* and *non-simple resonances*. This contrasts with the same system without the frequency modulation where resonance can occur only when the center frequency of the modulation coincides with the natural frequency of the undamped resonator, i.e., only when $\omega_0 = 1$.

### 4.1. The case of simple resonance – bounded responses

Taking $t \to \infty$, and invoking trigonometric identities, as steady state is reached all integrals in expression (11) can be written in terms of the following complex integral, which can be evaluated explicitly:

$$\int_0^\infty J_n(\delta\tau)\, exp(\mathrm{j}\omega\tau)\, d\tau = \begin{cases} \dfrac{exp[j\, n\, sin^{-1}(\omega/\delta)]}{\sqrt{\delta^2-\omega^2}}, & |\omega| < \delta \\ \dfrac{-j^{n+1}}{\sqrt{\omega^2-\delta^2}}\left(\dfrac{\delta}{\omega+\sqrt{\omega^2-\delta^2}}\right)^n, & \omega > \delta \\ \dfrac{j^{n+1}}{\sqrt{\omega^2-\delta^2}}\left(\dfrac{\delta}{-\omega+\sqrt{\omega^2-\delta^2}}\right)^n, & \omega < -\delta \end{cases} \quad (12)$$

The only non-decaying term in (12) corresponds to $\omega = 0$, and in terms of the solution (11) this translates to the relation $\omega_0 + \varepsilon n = \pm 1$, which recovers the resonance condition derived in the previous section for the damped case; provided that $n$ is integer, this case is referred to as *simple resonance*. The resulting steady-state response of (10) consists of quasi-periodic oscillations involving just two integer harmonics defined by, $n = (\pm 1 - \omega_0)/\varepsilon$, which, however, are activated (i.e., appear) at different time instants (as discussed previously, the $n$ −th harmonic is activated approximately at $\tau_n \sim n/\varepsilon\alpha$). Note that this analysis predicts this type of steady-state response for resonant ($\omega_0 = 1$), as well as off-resonant ($\omega_0 \neq 1$) center frequencies of the non-standard modulation. As an example, for the case of resonant center frequency, $\omega_0 = 1$, a simple asymptotic analysis based on (11) and (12) predicts the following steady-state response of the forced undamped oscillator (10):

$$x(t) \sim \underbrace{\frac{sint}{2a\varepsilon}\left[U(t) - U\left(t - \tau_{2/\varepsilon}\right)\right]}_{Harmonic\ n=0} + \underbrace{\frac{sint}{a\varepsilon}\, U\left(t - \tau_{2/\varepsilon}\right)}_{Harmonic\ n=-2/\varepsilon}, \quad t \gg 1 \quad (13)$$

where $\tau_{2/\varepsilon} \sim 2/\varepsilon^2\alpha$ is the activation time of the $n = -2/\varepsilon$ integer harmonic, and $U(s)$ is the Heaviside function Similar asymptotic expressions for the steady-state responses can be derived for modulated excitations with off-resonant center frequencies.

In Fig. 4 numerical simulations of the undamped resonator with wNSFM are shown, both for resonant (Fig. 4a) and off-resonant (Figs. 4b,c) modulation center frequencies. Starting from the case of resonant modulation center frequency, the numerical simulations agree with the asymptotic estimate (13) for the amplitudes of the steady-state responses, as well as the activation times of the two integer harmonics $n = 0$ and $n = -2/\varepsilon$. It is interesting to note the realization of large-amplitude quasi-periodic steady-state responses even in cases of off-resonant modulation center frequencies, both for below ($\omega_0 = 0.5$ – Fig. 4b) and above ($\omega_0 = 2.0$ – Fig. 4c) the natural frequency equaling unity of the undamped resonator. Clearly, in the absence of the weak modulation the off-resonance steady-state responses of the same resonator are small. Lastly, one observes the strongly modulated steady-state responses of the undamped oscillator (contrary to the resonant case $\omega_0 = 1$ where the steady-state response is nearly periodic – Fig. 4a); asymptotic estimates for these quasi-periodic responses can be derived from relations (11) and (12).

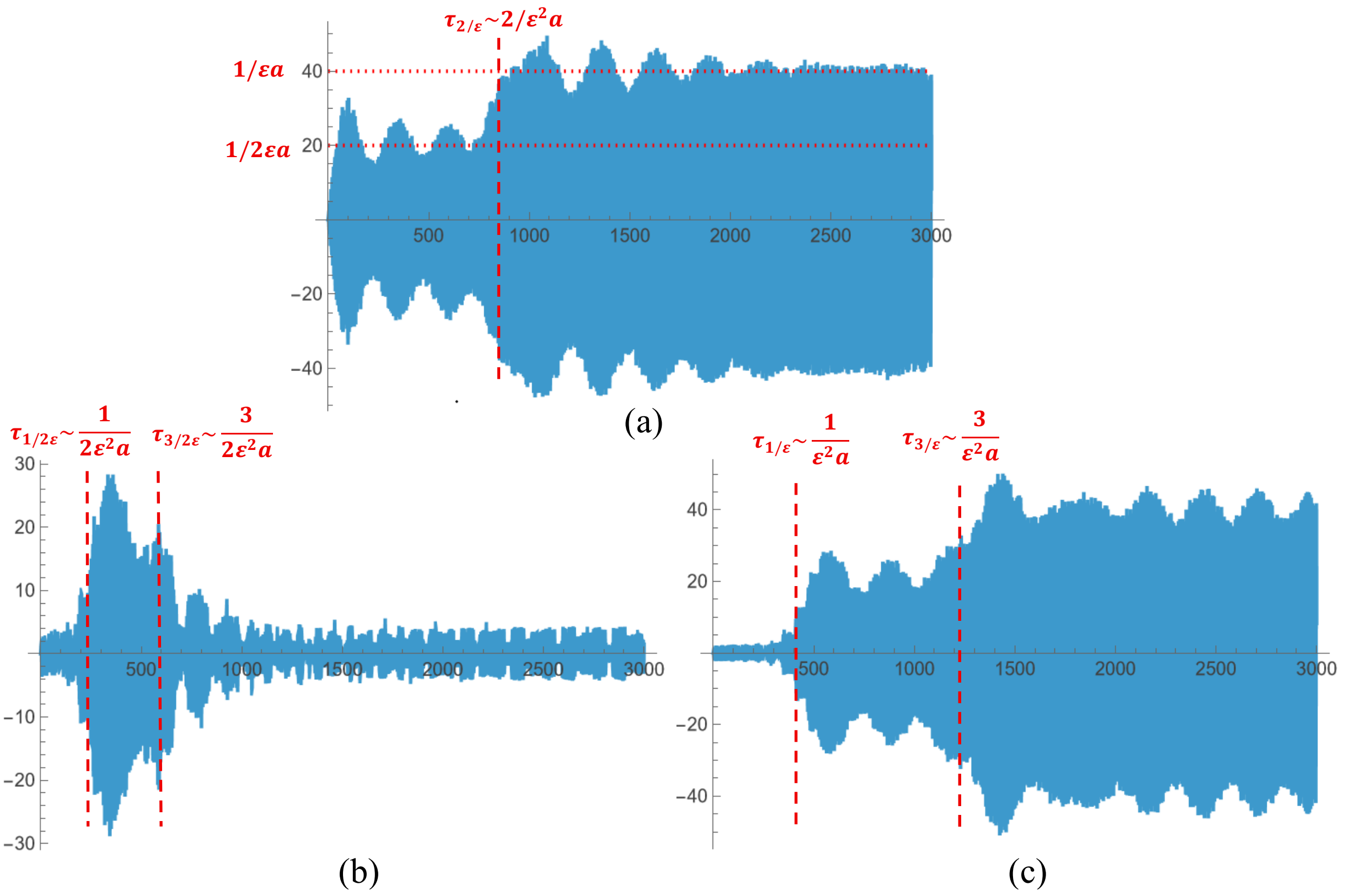

Figure 4. Steady-state responses of the undamped forced system (10) for the case of simple resonance, with $\varepsilon = 0.1$, $\alpha = 0.25$ and (a) resonant modulation center frequency $\omega_0 = 1$ (with the predicted asymptotic approximations for the steady-state amplitudes also shown), and off-resonant modulation center frequencies, (b) $\omega_0 = 0.5$, and (c) $\omega_0 = 2.0$.

### 4.2. The case of non-simple resonance – bounded or unbounded responses

When additional resonance conditions are met, the response of the undamped oscillator (10) may become unbounded as $t \to \infty$, resembling the classical growth to unbounded response of the same oscillator subject to unmodulated harmonic force with frequency equaling the natural frequency:

$$\ddot{x}(t) + x(t) = \text{cost}, \qquad x(0) = \dot{x}(0) = 0 \;\Rightarrow\; x(t) = (1/2)\, t\, sint \qquad (14)$$

To show this, the system (10) is reconsidered, and its long-time response is studied in more detail, making use of the long-term asymptotic approximation of the Bessel functions of first kind,

$$J_n(\varepsilon\alpha t) \sim \frac{1}{2}\sqrt{\frac{2}{\pi\varepsilon\alpha t}}\left[e^{j\left(\varepsilon\alpha t - \frac{n\pi}{2} - \frac{\pi}{4}\right)} + e^{-j\left(\varepsilon\alpha\tau t - \frac{n\pi}{2} - \frac{\pi}{4}\right)}\right], \qquad t \gg 1 \qquad (15)$$

so that the long-time asymptotic approximation of the modulated force in (10) is expressed as,

$$\ddot{x}(t) + x(t) \approx$$

$$\sum_{n=-\infty}^{\infty} \frac{1}{4}\sqrt{\frac{2}{\pi\varepsilon\alpha t}}\left[e^{j\left(\varepsilon\alpha t - \frac{n\pi}{2} - \frac{\pi}{4}\right)} + e^{-j\left(\varepsilon\alpha\tau t - \frac{n\pi}{2} - \frac{\pi}{4}\right)}\right]\left[e^{j(\omega_0+\varepsilon n)t} + e^{-j(\omega_0+\varepsilon n)t}\right] =$$

$$\sum_{n=-\infty}^{\infty} \frac{1}{4}\sqrt{\frac{2}{\pi\varepsilon\alpha t}}\left[e^{j\left[(\varepsilon\alpha-\omega_0-\varepsilon n)t - \frac{(2n+1)\pi}{4}\right]} + e^{j\left[(\varepsilon\alpha+\omega_0+\varepsilon n)t - \frac{(2n+1)\pi}{4}\right]} + cc\right], \tag{16}$$

$$t \gg 1$$

with initial conditions $x(0) = \dot{x}(0) = 0$, and $0 < \varepsilon \ll 1$. Hence, additional resonances in the undamped system can occur when the following conditions are satisfied:

$$\varepsilon\alpha \pm (\omega_0 + \varepsilon n) = \pm 1 \tag{17a}$$

These will be referred to as *non-simple resonances*, and the relations (17a) provide the *necessary* conditions for their realization. Additional *sufficient* conditions are now formulated.

The resonance conditions (17a) are further expressed as,

$$n = \pm\alpha - \frac{\omega_0 \pm 1}{\varepsilon} \tag{17b}$$

to highlight the fact that for non-simple resonance to occur the order of the harmonic, $n$, should be an integer. This is required since the Jacobi-Anger expansions (16) sum only over integer harmonics, restricting the domain of the order of Bessel functions and thus their frequency in the long-time approximations. Moreover, from (17b) for a given selection of the modulation parameters $\omega_0, \alpha$ and $\varepsilon$, at most four integer resonant harmonics may occur, although, as discussed below *non-simple resonance can occur with one, two or four integer resonant harmonics*, but not three integer resonant harmonics. To show this one denotes the indices of the four possible resonant harmonics as follows,

$$m = -\alpha - \frac{\omega_0 - 1}{\varepsilon}, \quad k = \alpha - \frac{\omega_0 + 1}{\varepsilon}, \quad p = -\alpha - \frac{\omega_0 + 1}{\varepsilon}, \quad q = \alpha - \frac{\omega_0 - 1}{\varepsilon}$$

without necessarily assuming that $m, k, p, q$ are integers. From these relations, it can be shown that the linear superposition of the (potentially) resonant harmonics must satisfy,

$$m + k - p - q = 0$$

It follows that, if one assumes that three of these resonant harmonics are integers, say $m, k, p \in \mathbb{Z}$, from the above relation the remaining resonant harmonic $q$ should be equal to a linear combination of integers $m, k, p$, so $q$ should be also integer, $q \in \mathbb{Z}$. Thus, the resonant response of the forced undamped oscillator may not contain only three resonant harmonics.

As the results in Fig. 5 show, the selection of the parameters of the frequency modulation dictates the existence of non-simple divergent resonance. Indeed, for $\omega_0 = 1, \alpha = 2, \varepsilon = 0.1$ (Fig. 5a) four resonant harmonics exist, $n = \{-22, -18, -2, 2\}$ and the undamped resonant response grows unboundedly. However, if the frequency modulation parameters are altered to prevent the formation of integer resonance resonant harmonics, $\omega_0 = 1, \alpha = 2.5$, $\varepsilon = 0.1$ (Fig. 5b) yielding the non-integers, $n = \{-22.5, -17.5, -2.5, 0.5\}$, the resulting undamped response is bounded.

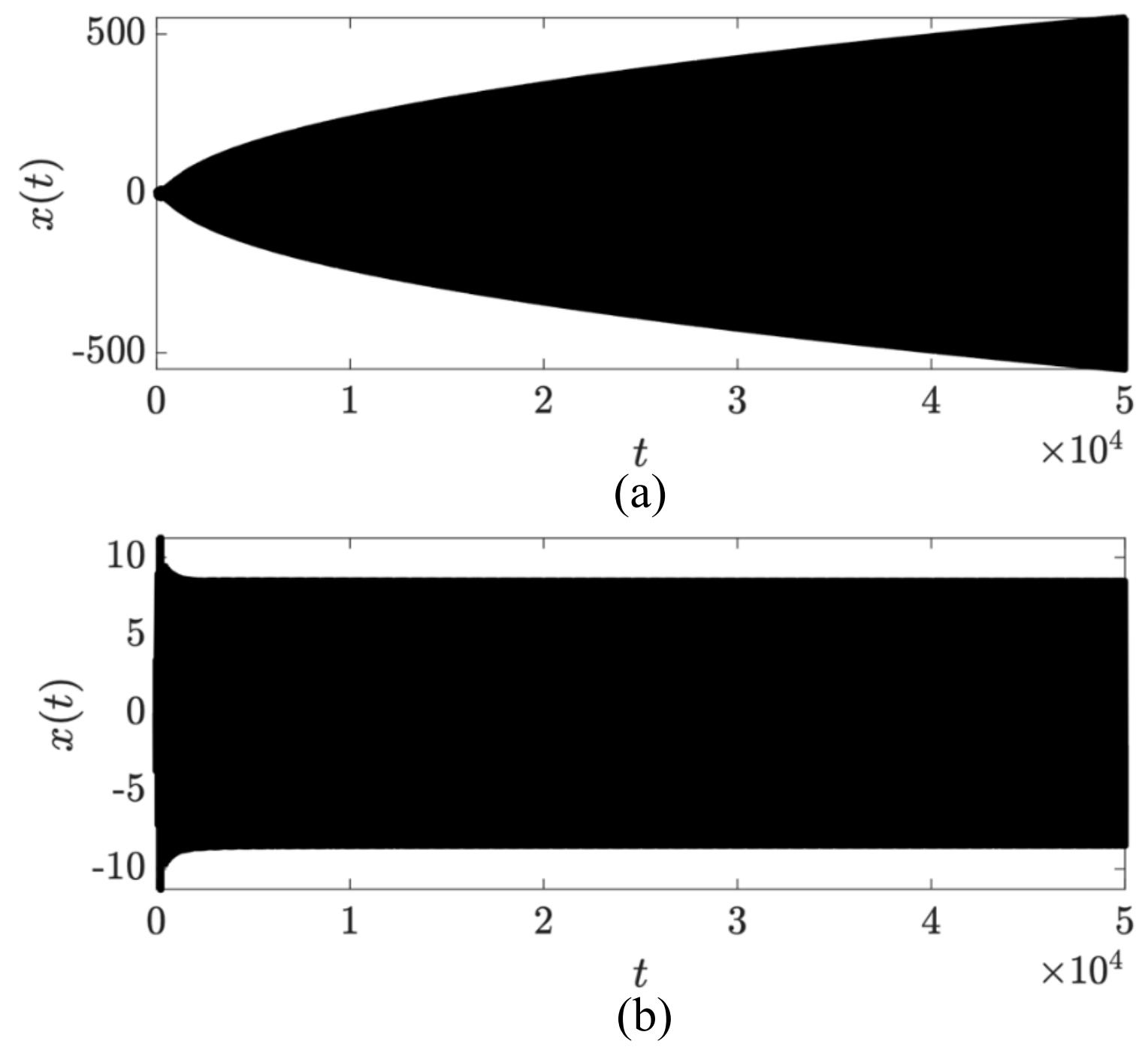

Figure 5. Responses of the undamped oscillator (10): (a) Divergent resonant response for $\varepsilon = 0.1, \omega_0 = 1, \alpha = 2$, (b) bounded non-resonant steady-state response for $\varepsilon = 0.1, \omega_0 = 1, \alpha = 2.5$.

From these results, one concludes that the non-simple resonance conditions (17a,b) alone, while necessary, are not sufficient to guarantee divergent (unbounded) growth in the response of the undamped oscillator (10). To formulate the necessary conditions, one reconsiders the asymptotic long-term approximation of the modulated forcing in (16) and notices that the oscillatory forcing terms have associated phases, $\phi_n = \pm(2n+1)\pi/4$. This suggests that the resonant harmonics in the modulated force can *destructively interfere* with each other, so that the resultant forcing terms will no longer cause the unbounded resonant excitations of the undamped oscillator. To prevent such destructive interference, the phases of the resonant harmonic should not be cancelled with each other, which provides the following additional condition for *constructive interference* of harmonics, and, hence, divergent (unbounded) resonance response of the undamped oscillator (10),

$$\sum_{n \in M} e^{j\phi_n} \neq 0, \tag{18}$$

where $M$ represents the subset of integers that correspond to the resonant harmonics, and $M = 1$, 2 or 4. This condition ensures that the resonant harmonics do not cancel out each other, and the resulting undamped resonant response has divergent growth as $t \gg 1$.

Summarizing the results of the previous analysis, the following two theorems are formulated.

<u>Theorem 1 (non-simple divergent resonance):</u> The necessary and sufficient conditions for the undamped oscillator (10) forced with weak non-standard frequency modulation to possess divergent resonant response as $t \gg 1$, are that the resonance conditions (17a,b) are satisfied for one, two or four integer harmonics $n$, and the constructive interference condition (18) is satisfied for the corresponding phases of these resonant harmonics. The resulting amplitude of the resonant undamped response grows as $\sim\sqrt{t}\,\sin\,(t+\mu)$ as $t \gg 1$.

Theorem 2 (non-simple non-divergent resonance): The necessary and sufficient conditions for the undamped oscillator (10) forced with weak non-standard frequency modulation to possess non-divergent (bounded) resonant response as $t \gg 1$, are that the resonance conditions (17a,b) are satisfied for two or four integer harmonics $n$, and the destructive interference, phase closure condition,

$$\sum_{n\in M} e^{j\phi_n} = 0 \quad \text{(Phase closure condition)}$$

is satisfied by the corresponding phases of the resonant harmonics, where $M$ denotes the set of resonant harmonics, and $M$ =2 or 4.

Some remarks are appropriate at this point. The last statement in Theorem 1 follows directly from the (slow) time dependence of the (common) amplitude of the oscillatory terms in the long-term approximation of the modulated forcing term in (16). Also, it is clear that only two or four resonant harmonics can produce the phase closure condition in Theorem 2, so the existence of a single resonant harmonic implies divergent resonance. In addition, the ratios of the frequencies of the resonant harmonics in Theorem 2 dictate the type of the bounded steady-state resonant response, and, generically, this ratio is expected to be irrational (yielding quasi-periodic long-term responses), rather than rational (yielding periodic steady-state responses).

In Fig. 6 three examples of divergent non-simple resonant responses are presented, corresponding to one (Fig. 6a), two (Fig. 6b) and four (Fig. 6c) resonant harmonics in constructive interference. For each simulation the corresponding vectorial phase summation (18) is also shown, indicating a non-zero resultant phase vector characterizing the divergent response; this results in non-zero average power input, $P_{\text{avg}}(t) = \frac{1}{t}\int_0^t \cos[(\omega_0 + \varepsilon\alpha\,\sin(\varepsilon\tau))\,\tau]\,\dot{x}(\tau)d\tau$, which appears as a non-smooth succession of "bursts" at every half cycle of the slow period, $T_\varepsilon = \pi/\varepsilon$, of the modulated force (see Fig. 6). These "bursts" manifest as trains of pulses with changing phases in in the corresponding plots of the instantaneous power input, $P_{in}(t) = \cos[(\omega_0 + \varepsilon\alpha\,\sin(\varepsilon t))\,t]\,\dot{x}(t)$.

Lastly, in Fig. 7 two examples of non-divergent (bounded) non-simple resonant responses are presented, corresponding to central frequencies $\omega_0 = 1$ (Fig. 7a), and $\omega_0 = 2$ (Fig. 7b) and four resonant harmonics in destructive interference. For each simulation the corresponding vectorial phase summation showing phase closure is also shown, yielding zero resultant phase that characterizes the bounded response as $t \gg 1$. In this case the corresponding average power input by the modulated force tends to zero with increasing time, $P_{avg}(t) \to 0$ as $t \gg 1$.

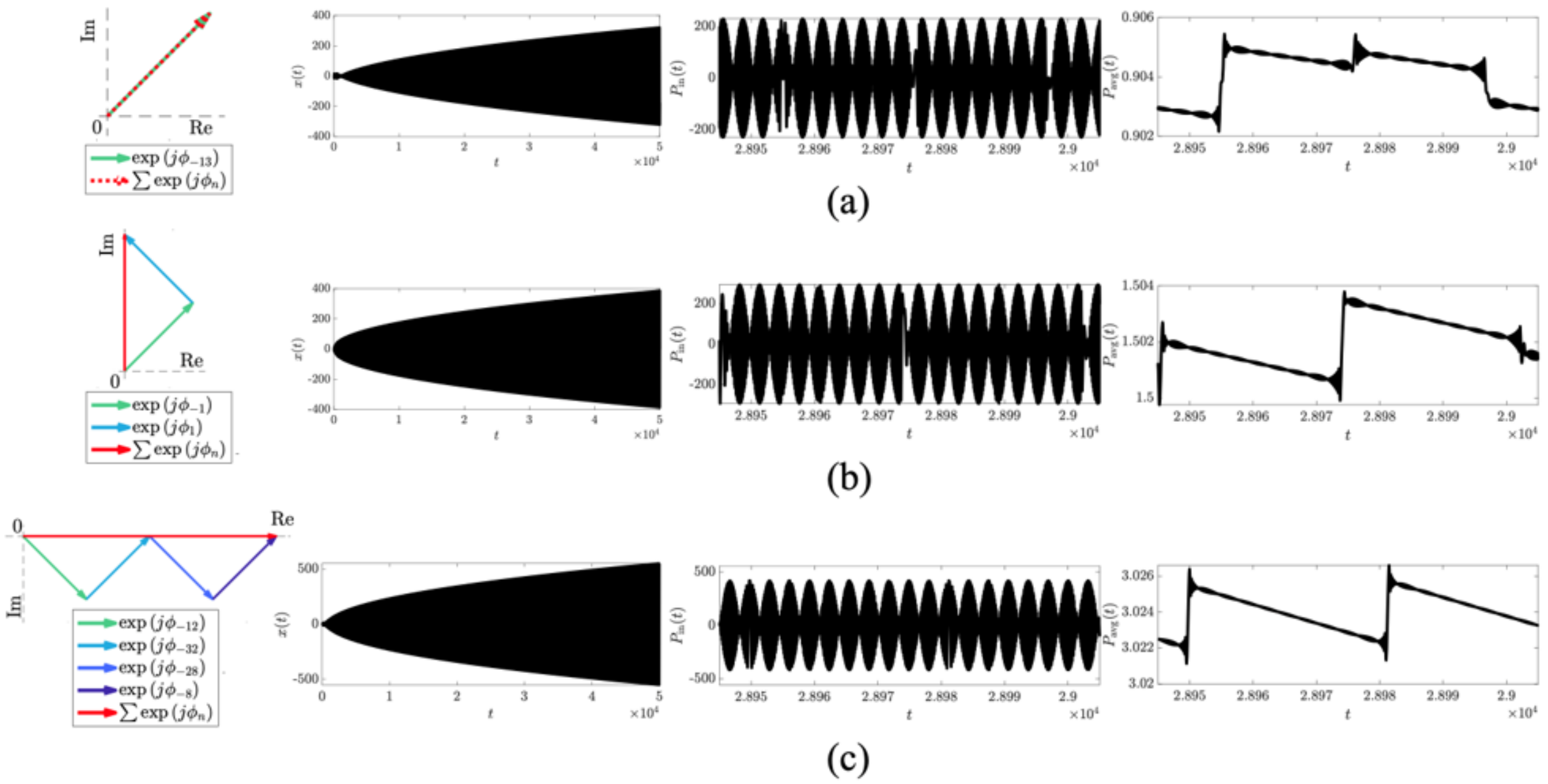

Figure 6. Non-simple resonances yielding divergent responses as $t \gg 1$: (a) $\varepsilon = 0.15, \alpha = 1/3, \omega_0 = 1$, single resonant harmonic $n = -13$, (b) $\varepsilon = 0.11, \alpha = 1, \omega_0 = 1$, two resonant harmonics, $n = -1, 1$, and (c) $\varepsilon = 0.1, \alpha = 2, \omega_0 = 2$, four resonant harmonics $n = -32, -28, -12, -8$; in each case the phase sum (18), the response, instantaneous power input, and the time-averaged power input are shown.

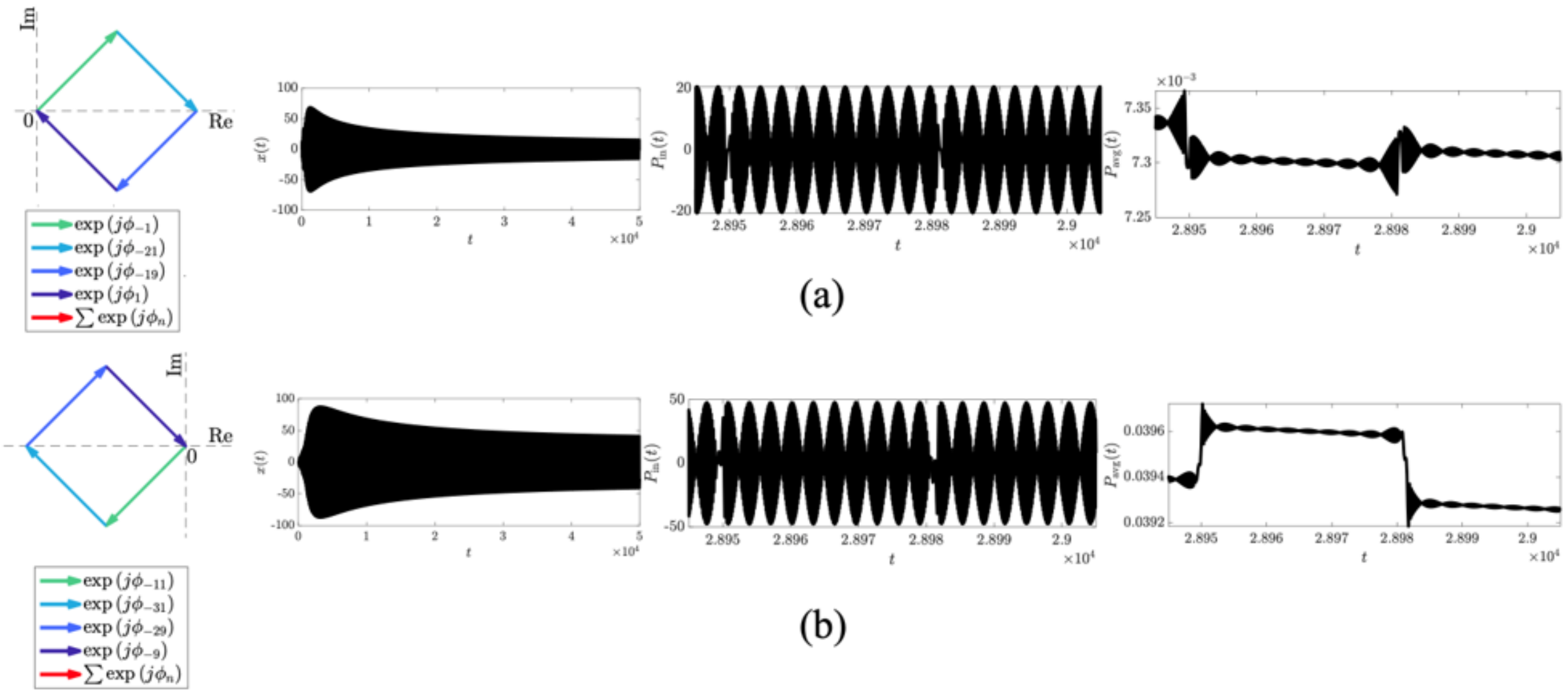

Figure 7. Non-simple resonances yielding non-divergent (bounded) responses as $t \gg 1$: (a) $\varepsilon = 0.1, \alpha = 1, \omega_0 = 1$, four resonant harmonics $n = -21, -19, -1, 1$, and (b) $\varepsilon = 0.1, \alpha = 1, \omega_0 = 2$, four resonant harmonics $n = -31, -29, -11, -9$; in each case the phase closure condition, the response, the instantaneous power input, and the time-averaged power input are shown.

In Table 1 below, we provide a synopsis of the excitation cases considered, and of the resulting long-term responses of the undamped oscillator (10) subject to wNSFM forcing. We note that the response of the classical undamped oscillator under harmonic (unmodulated) forcing does not fall into any of these classifications since the long-terms dynamics under the frequency modulation are different; for example, for the case of divergent non-simple resonance, the response of the undamped oscillator experiences unbounded growth with envelope $\sim\sqrt{t}$, compared to the similar envelope growth of the undamped oscillator under unmodulated harmonic forcing, which is $\sim t$.

**Table 1. Resonance in the undamped oscillator (10) with wNSFM forcing**

| Excitation Case | Resonance Conditions | Long-term Response |
|---|---|---|
| Simple Resonance | $\omega_0 + \varepsilon n = \pm 1$, $n \in \mathbb{Z}$ | $\sim \frac{1}{\alpha\varepsilon}$ (Bounded) |
| Non-simple Resonance – Non-divergent | $\varepsilon\alpha \pm (\omega_0 + \varepsilon n) = \pm 1$, $\sum_{n \in M} e^{j\phi_n} = 0$, $n \in \mathbb{Z}$ | Bounded |
| Non-simple Resonance – Divergent | $\varepsilon\alpha \pm (\omega_0 + \varepsilon n) = \pm 1$, $\sum_{n \in M} e^{j\phi_n} \neq 0$, $n \in \mathbb{Z}$ | $\sim\sqrt{t}$ (Unbounded) |

A last remark concerns the persistence of the studied non-simple divergent and non-divergent resonances, in the sense that they can be realized even for forced frequency modulations with central frequencies $\omega_0$ close or well apart from the natural frequency of the resonator (10). This indicates that by employing harmonic excitations with weak non-standard frequency modulations one can achieve a new type of resonances that persist to variations of the parameters of the applied harmonic forces. This is a unique new feature for forced linear time-invariant single-degree-of-freedom resonators, which, according to classical vibration theory can only be in resonance in frequency ranges close to their natural frequencies.

## 5. Concluding Remarks

Perhaps, the most interesting aspect of this work concerns the infinite set of resonances in the classical linear SDOF oscillator when excited by a harmonic force with the wNSFM (2). Contrary to classical single-frequency (unmodulated) harmonic forcing and traditional modulation patterns with constant frequency content – e.g., amplitude modulation or phase modulation (1), the forced response of the damped oscillator (5) always decays with increasing time (albeit very slowly), whereas the respective undamped oscillator can exhibit non-traditional divergence in its response due to non-simple resonance, with the envelope of its oscillatory response growing as square root of time. This non-simple resonance regime is persistent, in the sense that it is governed by the multi-parameter relationship (17a); that is, although a specific resonance can be eliminated, e.g., by varying the central frequency of the modulation, another resonance exists in the neighborhood of that resonance. Additional control of response divergence can be achieved through the phase relationships (18), allowing one to switch between divergent and non-divergent resonant responses by tuning the modulation parameters. It is indeed surprising that a signal of relatively simple functional form (2) can lead to such unexpected dynamical responses and to persistent infinite sets of resonances even in the simplest of linear oscillators. Indeed, the resonance relationships listed in Table 1 can be interpreted as a sort of combination resonances, which are well-known in nonlinear systems, but not in linear ones. In the current setting, the basic system is linear, but one observes resonant combinations of harmonics supplied by the external frequency-modulated forcing.

Physically, the wNSFM (2) generates a countable infinity of harmonics, each activated at a different time (with lower-frequency harmonics activated earlier). It follows that when a subset of such harmonics with appropriate phases and over the appropriate time windows are superimposed, their combination can produce a resultant harmonic with frequency component at or close to the natural frequency of the oscillator and thus give rise to resonance. In the damped system, realization of such a resonance occurs over a finite time window, resulting in a combination resonance capture. In the undamped system, constructive or destructive interference of a subset of harmonics (with appropriate phases) may result in bounded quasi-periodic or unbounded resonant responses.

In fact, many other functional forms of the frequency-modulated signals can provide time-variant frequency spectral content, and it is interesting to note that common (typically used) constant-frequency modulated signals should be considered as an exception. Therefore, one can hope to find even more complicated response regimes, resonant or non-resonant, allowing even broader control of the dynamics even for the simplest forced systems, like the class of linear oscillators considered herein.

The results presented in this work reveal a few interesting research venues moving forward. Extending the analysis to MDOF linear oscillators should be straightforward through standard modal analysis, whereby the forced response would be expressed as a superposition of SDOF modal responses for which the methodology developed here would apply. At this point we would like to conjecture that when a MDOF system is excited by a harmonic force with wNSFM having a central frequency close, say, to its first mode, it should be possible to induce delayed resonance excitation of a set of higher modes (depending on the bandwidth and the time-window of the applied frequency-modulated force); this should make it possible to replace standard frequency sweeps in modal testing with a simple harmonic excitation with wNSFM, an intriguing prospect indeed.

Extending this work to nonlinear oscillators would be more challenging, but the use of asymptotic methods, e.g., the method of multiple-scales, should enable a leading-order nonlinear study extending the resonance phenomena reported here for the linear oscillator (5); one of the anticipated challenges in the nonlinear case would be the frequency-energy relationship which would add an additional element of complexity to the analysis of the anticipated nonlinear resonance capture phenomena.

Regarding physical realizability, Chowning **[8]** addresses this issue in the context of accurately synthesizing the acoustics of musical instruments, requiring not only numerical generation of the input signal, digital-to-analog conversion, filtering, optimization through dynamic control of the time-varying modulation index to achieve desirable properties, and, finally, output to an audio amplifier and speakers. The present work, however, is focused solely on the dynamics of the system, has no optimization component, and otherwise follows the same analytical sequence, except that force to the vibrating system is applied by, e.g., an electromagnetic shaker. In the context of experimental implementation, a careful balance between the selection of the parameters of the frequency modulation and the practical implementation of the harmonic forcing with wNSFM needs to be considered, in terms, e.g., of the time-windows and the frequency sampling of the experimental measurements. However, we don't foresee any obvious limitations in terms of generating harmonic signals with wNSFM experimentally through standard shakers and frequency generators. However, in a practical setting the envelope of the phase should be restricted below a definite bound, which, in turn, would restrict the bandwidth of the modulated force. Given, however, the analytical findings reported in this work, one could predictively select the parameters of the wNSFM to target specific frequency ranges and time windows for an experiment. So, the proposed new type of frequency modulation introduced in this work should be practically realizable in, e.g., modal testing applications.

## Acknowledgments

The authors would like to express their appreciation to the reviewers of this work for their thorough reviews and insightful comments which helped significantly the authors revise the manuscript, relating this work to previous works, and correct some errors in the original submission. O.V.G is grateful to the Israel Science Foundation (grant 3085/25) for financial support.